\documentclass[a4paper,12pt]{article}
\usepackage{amsmath, amssymb}
\usepackage{doc, exscale, fontenc, latexsym, syntonly}
\usepackage{xypic}

\author{ A. Mutlu, T. Porter}
\title{\large{\bf Free crossed resolutions from simplicial resolutions with given $CW$-basis}} 
 
\newtheorem{defn}{Definition}[section]
\newtheorem{prop}[defn]{Proposition}
\newtheorem{thm}[defn]{Theorem}
\newtheorem{lem}[defn]{Lemma}
\newtheorem{cor}[defn]{Corollary}
\newenvironment{pf}{{\bf Proof:}}{\hfill$\Box$\mbox{}}

\date{}

\begin{document}


\maketitle

\begin{flushright}
\parbox{11.3cm}{
\centerline{\small\bf R\'esum\'e}
{\small
Dans cet article, les Auteurs examinent la relation entre une $CW$-base pour un
groupe simplicial, et des m\'ethodes \`a engendrer librement le complexe crois\'e
associ\'e. On examine en d\'etail le cas des r\'esolutions, en comparant les
r\'esolutions simpliciales libres et les r\'esolutions crois\'ees d'un groupe.

{\noindent\bf  Classification A. M. S. :}  18D35, 18G30, 18G50, 18G55, 20F05,
57M05.}}
\end{flushright}

\section*{Introduction}

When J.H.C. Whitehead wrote his famous papers on ``Combinatorial
Homotopy'', \cite{jh1}, it would seem that his aim was to produce
a combinatorial, and thus potentially constructive and
computational, approach to homotopy theory, analogous to the
combinatorial group theory developed earlier by Reidemeister and
others.  In those papers, he introduced CW-complexes and also the
algebraic `gadgets' he called homotopy systems, and which are now
more often called free crossed complexes, \cite{bh1}, or totally free crossed
chain complexes, \cite{baues1}.

Another algebraic model for a (connected) homotopy type is a simplicial group 
and again, there, one finds a notion of freeness.  In both cases we have 
`freeness', yet no easily defined category of things on which our objects
are `free'.  Kan, \cite{kan2}, introduced the notion of a CW-basis for a free simplicial
group and more recently, \cite{rab}, R.A.Brown has introduced Peiffer-Whitehead
word systems or extended group presentations as a means of presenting a 
`homotopy system'.  In both cases the aim is to use the `generators' as a 
combinatorial way of controlling or manipulating the algebraic model, i.e. 
extending the `yoga' of combinatorial group theory to higher dimensions.

There are ways of passing from a simplicial group to a crossed complex 
(see for example, \cite{ep}) and as these are all equivalent to a left 
adjoint, one expects freeness to be preserved, and it is, but this is not 
trivial.  As we do not know on what type of thing the simplicial group is 
free, nor on what the crossed complex is free, the conclusion is not a
simple consequence of left adjointness of some sort.  The problem is
that to construct the $n$\textsuperscript{th} level, you need some
generators together with a map to the $(n-1)$\textsuperscript{st} one, and
of course you cannot do that until that level is constructed!

In this paper we apply methods from our earlier papers
\cite{mp,mp1}, to examine the relationship between the notions of
free basis for simplicial groups and that for crossed complexes.  We
have included a shortened proof of the result from \cite{carrasco} and
\cite{ep}, describing the passage from simplicial groups to crossed
complexes, as this allows for a direct verification of freeness at the
base of the crossed complex.

Although our results would seem to apply in general, we have restricted
detailed attention to simplicial resolutions.  This is partially since
there are known problems of non-realizability of a homotopy system by a
CW-complex (cf. Whitehead, \cite{jh1}, section 15, or R.A.Brown,
\cite{rab}, p.527) and hence by a free simplicial group with CW-basis.
It thus seems prudent to understand these non-realizability  results
better from this simplicial viewpoint before attempting to look at the
general case.  Those results do not seem to disturb the general case in
any significant way, but they leave them somewhat incomplete in the view
they give of the general problem.

\section{Preliminaries }

We will denote the category of  groups by $\mathfrak{Grp}.$


\subsection{Simplicial Groups}


A simplicial group ${\bf G}$ is a simplicial object in the category of groups.
We will denote the category of simplicial groups by $\mathfrak{SimpGrp}.$
We will only need a small amount of the extensive theory of simplicial groups 
here and would refer to the book by May \cite{may} or the survey by Curtis \cite{curtis}
for information on the more `classical' parts of the theory. We will assume a basic 
knowledge of the elementary homotopy of simplicial sets and simplicial groups, 
but will also refer to facts and concepts from earlier parts of this series of 
papers, \cite{mp,mp1,mp2}.  If ${\bf G}$ is a simplicial group, then $(NG,\partial)$ will be the corresponding Moore complex.  Our conventions on this
and related notions are given in \cite{mp}.


\subsection{Step By Step Constructions}


This section is a brief r\'esum\'e of how to construct simplicial
resolutions. The work depends heavily on a variety of sources, mainly 
\cite{andre}, \cite{keune} and \cite{moore}. Andr\'e only treats commutative
algebras in detail, but Keune \cite{keune} does discuss the general case quite 
clearly.


First recall the following notation and terminology which will be used in
the construction of a simplicial resolution.

Let $[n]$ be the ordered set, $[n]=\{0<1<\cdots<n\}$. We define the following
maps: Firstly the injective monotone map $\delta _i^n:[n-1]\rightarrow [n]$
is given by 
$$
\delta _i^n(x)=\left\{ 
\begin{array}{lcl}
x & \text{if} & x<i, \\ 
x+1 & \text{if} & x\geq i, 
\end{array}
\right. 
$$
for $0\leq i\leq n\neq 0.$
An increasing surjective monotone map 
$\alpha_{i}^{n} : [n+1] \to [n]$ is  given by 
\[ \alpha_{i}^{n}(x) = \left\{ \begin{array}{ll} 
                         x    &  \mbox{ if $ x \leq i, $ } \\
                         x-1  &  \mbox{ if $ x > i, $ }
                        \end{array}
                        \right. \]
for $0 \leq i \leq n$. 
We denote by $\{m,n\}$ the set of increasing surjective maps $[m]\rightarrow
[n].$\\ 

\medskip

\textbf{Killing Elements in Homotopy Groups}\\ 
The following section describes the `step-by-step' construction due to Andr\'e 
\cite{andre}, that source however concentrates on simplicial algebras. We have 
adapted his treatment to handle simplicial  groups.

We recall that if $F$ and $G$ are groups, a map 
\[ G \times F \to F \]
\[(g,f) \longmapsto{}^{g} f \]
is a left action if and only if for all $g,{g'} \in G, f,{f'}\in G,$

$1.$  ${}^{g}(f {f'}) ={}^{g}f{~}{}^{g}{f'},$

$2.$  ${}^{g{g'}}f = {}^{{g}}({}^{g'}f),$

$3.$  ${}^{1}f = f.$ \\
In this case we say $F$ is a $G$-group.

Let {\bf G} be a simplicial group and let $k\geq 1$ be fixed. Suppose we are given a set $\Omega $ of elements $
\Omega = \{x_\lambda :\lambda \in \Lambda \}$, $x_\lambda \in \pi _{k-1}({\bf G}),$ then we can choose a corresponding set
of elements ${\it \vartheta}_\lambda \in NG_{k-1}$ so that
$x_\lambda ={\it \vartheta}_\lambda{~}\partial _k(NG_k).$
(If $k=1,$ then as $NG_0=G_0,$ the condition that ${\it \vartheta}_\lambda \in NG_0$
is immediate.) We want to `kill' the elements in $ \Omega$.

We form a new simplicial group $F_n$ where

1) $F_n$ is the free $G_n$-group, 
$$
F_n= \coprod_{\lambda,t}G_n\{y_{\lambda ,t}\}\ \text{with}\ \lambda \in \Lambda \text{ and }\ t\in
\{n,k\}, 
$$
where $G_n\{y\} =G_n*<y>,$ the free product of $G_n$ and a free group generated by $y.$

2) For $0\leq i\leq n$, the group homomorphism $s_i^n:F_n\rightarrow
F_{n+1}$ is obtained from the homomorphism $s_i^n:G_n\rightarrow G_{n+1}$
with the relations

\begin{center}
$s_i^n(y_{\lambda ,t})=y_{\lambda ,u}$ \ \ with \ \ $u=t\alpha_i^n,\ \
t:[n]\rightarrow [k].$
\end{center}

3) For $0\leq i\leq n\neq 0,$ the group homomorphism $d_i^n:F_n\rightarrow
F_{n-1}$ is obtained from $d_i^n:G_n\rightarrow G_{n-1}$  with the
relations 
$$
d_i^n(y_{\lambda ,t})=\left\{ 
\begin{array}{clcl}
y_{\lambda ,u} & \text{if the map} & u=t\delta _i^n & \text{is surjective},
\\ t{^{\prime }}({\it \vartheta}_\lambda ) & \text{if } & u=\delta _k^kt{^{\prime }},
&  \\ 
1 & \text{if } & u=\delta _j^kt{^{\prime }} & \text{with}\ \ j\neq k, 
\end{array}
\right. 
$$
by extending multiplicatively. 

We sometimes denote the ${\bf F} $ so constructed by ${\bf G}(\Omega)$.

\medskip

{\bf Remark : } 
In a `step-by-step' construction of a simplicial resolution,
(see below), there are thus
the following properties: i) $F_n= G_n$ for \mbox{$n<k$,}
ii) $F_k=$  a free $G_k$-group over a set of
non-degenerate indeterminates, all of whose faces are the identity except the
$k^{th}$, and
iii) $F_n$ is a free $G_n$-group on some degenerate
elements for $n>k.$

We have immediately the following result, as expected.

\begin{prop}\label{veli}
The inclusion of simplicial groups ${\bf G}\hookrightarrow{\bf F}$,
where ${\bf F}={\bf G}(\Omega )$, induces a homomorphism 
$$
\pi_n({\bf G})\longrightarrow \pi_n({\bf F}) 
$$
for each $n$, which for $n<k-1$ is an isomorphism,
$$
\pi_n({\bf G})\cong \pi_n({\bf F}) 
$$
and for $n=k-1$,  is an epimorphism with kernel generated
by elements of the form $\bar \vartheta_\lambda = \vartheta_\lambda \partial _kNG_k$,
where $\Omega = \{x_{\lambda}: \lambda\in\Lambda\}.$
\end{prop}\hfill$\square$

\textbf{Constructing Simplicial Resolutions}\\
The following result is essentially  due to Andr\'{e} \cite{andre}.
\begin{thm}
If $G$ is a  group, then it has a free simplicial resolution 
$\mathbb{F}$.
\end{thm}
\begin{pf}
The repetition of the above construction will give us the simplicial resolution
of a group. Although `well known', we sketch the construction so as to 
establish some notation and terminology.

Let $G$ be a  group. The zero step of the construction consists of a
choice of a free group F and a surjection $g : F\rightarrow G$
which gives an isomorphism $F/{\rm Ker}g\cong G$ as groups. Then
we form the constant simplicial group ${\mathbf  F}^{(0)}$ for which in every
degree $n,$ $F_n=F$ and $d_i^n=$ id $=s_j^n$ for all $i,j.$ Thus 
${\mathbf F}^{(0)}={\bf K}(F,0)$ and $\pi _0({\mathbf F}^{(0)})=F.$
Now choose a set $\Omega^0$ of normal generators of the normal subgroup 
$N={\rm Ker}(F\stackrel{g}{\longrightarrow}G),$ and obtain the simplicial group in 
which $F_1^{(1)}=F(\Omega ^0)$ and for $n>1,$ $F_n^{(1)}$ is a free 
$F_n$-group over the degenerate elements as above. This simplicial group will be denoted by 
$\bf{F}^{(1)}$ and will be called the {\em 1-skeleton of a simplicial
resolution of the group} $G$.

The subsequent steps depend on the choice of sets, $\Omega^0$, 
$\Omega^1, \Omega^2, \ldots, \Omega^k,\ldots .$ Let ${\mathbf F}^{(k)}$ be the
simplicial group constructed after $k$ steps, the $k$-skeleton of the
resolution. The set $\Omega^k$ is formed by elements $a$ of $F_k^{(k)}$
with $d_i^k(a)=1$ for $0\leq i\leq k$ and whose images $\bar a$ in 
$\pi _k({\mathbf F}^{(k)})$ generate that module over $F_k^{(k)}$ and ${\bf F}^{(k+1)}$.

Finally we have inclusions of simplicial groups 
$$
{\mathbf F}^{(0)}\subseteq {\mathbf F}^{(1)}\subseteq\cdots \subseteq 
{\mathbf F}^{(k-1)}\subseteq {\mathbf F}^{(k)}\subseteq \cdots  
$$
and in passing to the inductive limit (colimit), we obtain an acyclic free
simplicial group ${\mathbf F}$ with ${D}_n = F_n^{(k)}$ if $n\leq k.$
${\mathbb F} = ({\mathbf F},g)$ is thus a simplicial resolution of the
 group $G$. 

The proof of theorem is completed.~\end{pf}

\medskip

{\bf Remark : }
A variant of the `step-by-step' construction gives: if     
${\bf G}$ is a simplicial group, then there exists a free
simplicial group ${\mathbf F}$  and an epimorphism
$
{\mathbf F}\longrightarrow {\bf G} 
$
which induces isomorphisms on all homotopy groups. The details are
omitted as they are well known.

\medskip

{\bf Terminology : } It is sometimes useful to write ${\mathbb F}^{(k)} =
({\mathbf F}^{(k)},g)$ for the augmented simplicial group constructed at
the $k^{th}$ step. 
The data needed to go from $\mathbb{F}^{(k)}$ to $\mathbb{F}^{(k+1)}$ are more 
precisely a
set $\Omega^k$ and a function $g^{(k)}: \Omega^{k}\longrightarrow F_k^{(k)}$ whose image is
contained in $NF_k^{(k)}$ and which generates $\pi_k(\mathbb{F}^{(k)}).$~ (We often consider
$g^{(k)}$ as being an inclusion and leave it out of the notation.) The pair $(\Omega^k,
g^{(k)})$ is then called {\em k-dimensional construction data for the resolution} and the
finite sequence 
$
((\Omega^0,g^{(0)}),\ldots,(\Omega^{k-1},g^{(k-1)}))
$
is called a $k$\textsuperscript{th}-{\em level presentation} of the group $G$.

The key observation, which follows from the universal property of the construction, 
is a freeness statement:
\begin{prop}\label{free}
Let $\mathbf{F}^{(k)}$ be a $k$-skeleton of a simplicial resolution of $G$ and 
$(\Omega^k, g^{(k)})$ $k$-dimension construction data for $\mathbb{F}^{(k+1)}.$ 
Suppose given a simplicial group
morphism 
$\Theta:\mathbf{F}^{(k)}\longrightarrow {\bf G}$ such that $\Theta_{\ast}
(g^{(k)}) =0,$ then $\Theta$ extends over $\mathbb{F}^{(k+1)}.$ 
\end{prop}

This freeness statement does not contain a uniqueness clause. That can be achieved by
choosing a lift for $\Theta_kg^{(k)}$ to $NG_{k+1},$ a lift that must exist since
$\Theta_{\ast}(\pi_k(\mathbb{F}^{(k)}))$ is trivial.

When handling combinatorially defined resolutions, rather than functorially defined ones,
this proposition is as often as close to `left adjointness' as 
is possible without entering the realm of homotopical algebra to an extent greater than is
desirable for us here.

We have not talked here about the homotopy of simplicial group
morphisms, and so will not discuss homotopy invariance of this construction
for which one adapts the description given by Andr\'e, ~\cite{andre}, or Keune, \cite{keune}.

\subsection{CW-bases}

We recall from \cite{kan2} and \cite{curtis} the following definitions.
\begin{defn}
A simplicial group {\bf F} will be called free if\\
(a)\qquad $F_n$ is a free group with a given basis, for every integer $n\geq 0,$\\
(b)\qquad The bases are stable under all degeneracy operators, i.e., for every pair 
of integers $(i,n)$ with $0\leq i\leq n$ and every given generator $x\in F_n$ the 
element $s_i(x)$ is a  given  generator of $F_{n+1}.$
\end{defn}
\begin{defn}
Let ${\bf F}$ be a free simplicial group (as above). A subset $\mathfrak{F}\subset {\bf F}$ 
will be called a $CW-basis$ for ${\bf F}$ if \\
(a)\qquad $\mathfrak{F_n} = \mathfrak{F}\cap F_n$ freely generates 
$F_n$ for all $n\geq 0,$\\
(b)\qquad $\mathfrak{F}$ is closed under degeneracies, i.e., $x\in \mathfrak{F_n}$
implies $s_i(x)\in \mathfrak{F_{n+1}}$ for all $0\leq i\leq n,$\\
(c)\qquad if $x\in\mathfrak{F_n}$ is non-degenerate, then $d_i(x) = e_{n-1}$, 
$(e_{n-1}$, the identity element of $F_{n-1}$) for all $0\leq i< n.$
\end{defn}
Let $\bf{F}$ be a free simplicial group with given $CW$-basis, 
$\mathfrak{F},$ then $X_0=\mathfrak{F}_0$ freely generates
$F_0$, that is, $F_0 = F(X_0).$  In general, note that if $Y_n = \mathfrak{F}_n\setminus\bigcup\limits_{i=0}^{n-1}s_i(\mathfrak{F}_i)$
then $Y_n\subseteq NF_n$.

\subsection{Crossed Modules}
J. H. C. Whitehead $(1949)$ \cite{jh1} described crossed modules in various 
contexts especially in his investigation into the group structure of 
relative homotopy groups. 
\begin{defn}
A {\rm pre-crossed module of
 groups} consists of a group, $G_1,$ a $G_1$-group $G_2$, and a group homomorphism 
$
\partial :G_2\longrightarrow G_1, 
$
such that for all $g_2 \in G_2,g_1 \in G_1$ 
$\\
\begin{array}{cccc}
CM1)\ \ \ \  & \partial ({}^{g_1}g_2) & = & g_1\partial(g_2){g_1}^{-1}.
\end{array}
$\\
This is {\rm a crossed module }if in addition, for all $g_2,g_2^{\prime }\in G_2$,\\
$
\begin{array}{cccc}
CM2)\ \ \ \  & {}^{\partial(g_2)}{g_2'} & = & g_2{g_2'}(g_2)^{-1}.
\end{array}
$
\end{defn}
The second  condition (CM2) is called {\em the Peiffer identity}. We denote such a
crossed module by $(G_2,G_1,\partial )$. Clearly any crossed module is a
pre-crossed module.

\subsection{Free Crossed Modules}
The notion of a free crossed module was  
described by J. H. C. Whitehead ~\cite{jh1}. We refer the reader to \cite{bh3} 
for the construction of a free crossed module on a presentation and the proofs 
of the results below.  The related notion of totally free (pre-)crossed module 
is discussed in \cite{mp2}.

\begin{thm}\label{rin}
A free crossed module $G_1$-module $(G_2,G_1,\partial )$ exists on any function $f:S\rightarrow G_1$ with codomain $G_1.$
\end{thm}
\begin{pf}
See ~\cite{bh3}.~
\end{pf}

\medskip

If ($G_2,G_1,\partial $) is a free crossed $G_1$-module on the trivial function $$1:S\rightarrow G_1,$$
then $G_2$ is a free $G_1$-module on the set $S$.

\section{Crossed Complexes}

\subsection{Peiffer pairings and boundaries in the Moore complex}
Firstly we recall from ~\cite{mp} the following result.
Let
${\bf G}$ be a simplicial group with Moore complex ${\bf NG}$ and  for $n \geq 1$, let
$D_{n}$ be the normal subgroup generated by the degenerate elements of $n$. 
If $G_{n} \ne D_{n}$, then 
\[ NG_{n}\cap D_n = N_{n} \cap D_n \quad\mbox{ for all $n \geq 1,$ } \]
where $N_n$ is a normal subgroup in $G_n$ generated by an explicitly given 
set of elements.

\subsection{Crossed Complexes and Crossed Resolutions}
The definition of a crossed complex (over a groupoid) was  given by
R. Brown and P. J. Higgins (1981)~\cite{bh1} generalising earlier work
of Whitehead (1949)~\cite{jh1}.  Crossed resolutions are discussed in several
sources. A particularly useful one is the thesis of Tonks, \cite{tonks}, which
handles constructions of crossed resolutions in some detail.

 \subsection{Peiffer-Whitehead word systems}

R.A.Brown, \cite{rab}, introduces a system of generators for a `homotopy system' 
that he calls a {\em Peiffer-Whitehead word system}.  His sets of generators 
are only in a finite number of 
dimensions whilst ours may need to be in an infinite set of levels to get a 
resolution, so his needs are not the same as ours, but nonetheless it
seems worthwhile to include his definition as it provides a point of
comparison with his work:

\begin{defn}{\rm \cite{rab}, p. 525}

A {\em Peiffer-Whitehead word system} or {\em extended group presentation }
$ W$ consists of a finite list of finite sets  $\langle
W^{(1)}|W^{(2)}|\ldots W^{(n)}\rangle$ together with {\em boundary
homomorphisms} $d_3, \ldots, d_n$ described as follows:\\
$W^{(1)} = I_1$ is a set of indices;\\
$W^{(2)} = \{w^2_\beta|\beta \in I_2\}$ is a set of words representing
elements of the free group $F = F(I_2)$;\\
$W^{(3)} = \{w^3_\gamma|\gamma \in I_3\}$ is a set of words representing
elements of the free $F$-crossed module $C(I_2)$ with boundaries
$\{c_\beta =\langle w^2_\beta\rangle | \beta \in I_2\}$;\\
$W^{(m)} = \{w^m_\mu|\mu \in I_m\}\quad (4< m \leq n)$ is a set of words 
representing elements of the free $\mathbb{Z}G$-module  $M_m = M_m(I_{m-1})$, where
$G$ is the group presented by $\langle W^{(1)}|W^{(2)}\rangle$;\\
$d_3 : C(I_2) \rightarrow F$ is a group homomorphism determined by $d_3
(i_\beta) = \langle  w^2_\beta \rangle$;\\
$d_4 : M_4(I_3) \rightarrow C(I_2)$ is a homomorphism determined by $d_4
(i_\gamma) = \langle  w^3_\gamma \rangle$;\\
$d_m : M_m(I_{m-1}) \rightarrow  M_{m-1} \quad (4< m \leq n)$ is a module 
homomorphism determined by $d_m(i_\lambda) = \langle  w^{m-1}_\lambda \rangle$;\\
In addition all words must have trivial boundaries 
$$d_m\langle w^m_\mu \rangle = identity \mbox{ \quad for } \mu \in I_m,  3\leq m
\leq n.$$
\end{defn}
Such a word system clearly specifies the generators of each level and
their images in the next level down.

\subsection{From Simplicial Groups to Crossed Complexes}\label{fa}
P. Carrasco and A. M. Cegarra ~\cite{c:c} defined 
$$
C_n({\bf G})=\frac{NG_n}{(NG_n\cap D_n)~d_{n+1}(NG_{n+1}\cap D_{n+1})} 
$$
for a simplicial group {\bf G}.
This gives a crossed complex ${\mathfrak{C}}({\bf G})$ starting from the Moore complex ({\bf NG}, $
\partial )$ of ${\bf G}.$ The map 
$
\partial _n:C_n({\bf G})\rightarrow C_{n-1}({\bf G})  
$
will be that induced by $d_n^n$.
Their proof requires an understanding of hypercrossed
complexes. P. J. Ehlers and T. Porter, ~\cite{ep}, developed a more direct proof for
simplicial groupoids. Here we will sketch a shorter
argument showing that $ {\mathfrak{C}}({\bf G})$ is a
crossed complex as we can use  some of the ideas later on. This proof
emphasises the role played by the various $F_{\alpha,\beta}.$ These pairing operations on the 
Moore complex were introduced in \cite{mp} and \cite{mp1}.  They are defined
by forming $[s_\alpha x,s_\beta y]$ and then projecting the result into
the Moore complex.  Detailed examples and calculations are given in these
papers cited above.
If $x\in NG_n$, we will write $\bar x$ for the corresponding element of 
$C_n({\bf G})$. 
\begin{lem}
The subgroup $(NG_n\cap D_n)~ d_{n+1}(NG_{n+1}\cap D_{n+1})$ is a normal 
subgroup in 
$G_n$.
\end{lem}
\begin{pf}
This is a routine use of the degeneracies.
\end{pf}
\begin{prop}
Let {\bf G} be a simplicial group, then defining
$C({\bf G})= (C_n({\bf G}), \partial)$ as above
yields a crossed complex. 
\end{prop}
\begin{pf}
(i) For $n\geq 2,$ $C_n({\bf G})$ is abelian, in fact 
$$
\begin{array}{rcl}
F_{(n-1)(n)}(x,y)& = &[s_{n-1}x, ~s_ny]~[s_ny, ~s_nx]\\
d_{n+1}F_{(n-1)(n)}(x,y)& = & [x, ~s_{n-1}d_n(y)]~[y,x] \end{array}$$
is in $(NG_n\cap D_n)d_{n+1}(NG_{n+1}\cap D_{n+1}), $
so $d_{n+1}F_{(n-1)(n)}(x,y) \equiv  1 $ mod $(NG_n\cap D_n)d_{n+1}(NG_{n+1}\cap D_{n+1})$
giving $\overline{x}\overline{y} = \overline{y}\overline{x}.$\\ 
(ii) For $x\in NG_n,$ and $y\in NG_m,$  taking $\alpha=(n,n-1,\ldots,m)$, and $\beta =(m-1),$
it is easy to see that\\
$$
F_{\alpha, \beta}(x_{\alpha}, y_{\beta}) = \prod\limits_{k=0}^{n-m+1}
[s_ns_{n-1}\ldots s_mx, ~s_{m-1+k}y]^{(-1)^k} [s_ns_{n-1}\ldots s_m(x),~s_{n}y]$$
$$d_{n+1}F_{\alpha, \beta}(x_{\alpha}, y_{\beta})  = \prod\limits_{k=0}^{n-m} [s_ns_{n-1}\ldots s_mx, ~s_{m-1+k}y]^{(-1)^k} [s_{n-1}\ldots s_m(x),~y].
$$
This implies that
$[s_m^{(n-m)}(x), y] \in (NG_n\cap D_n)d_{n+1}(NG_{n+1}\cap D_{n+1}),$
(where $s_{m}^{(n-m)}x = s_{n-1}\ldots s_mx$) which shows that the actions of $NG_m$ on $NG_n$
defined by conjugation
$$
{}^{\overline{x}}\overline{y}= \overline{s_{m}^{(n-m)}(x)ys_{m}^{(n-m)}(x)^{-1}}
$$
via these degeneracies are trivial if $m\geq 1.$ For $m=1,$ this gives $\alpha =(n,n-1,\ldots,1)$, $\beta=(0)$
and\\
\hspace*{1.5cm}$
F_{(n,n-1\ldots,1)(0)}(x, y) = 
$\\
\hspace*{3cm}$
\prod\limits_{k=0}^{n}
[s_ns_{n-1}\ldots s_1x, ~s_{0+k}y]^{(-1)^k}
[s_ns_{n-1}\ldots s_1(x),~s_{n}y],$
\\where $x\in NG_1, y\in NG_n,$ and it is easily checked that\\ 
\hspace*{1.5cm}
$
d_{n+1}F_{(n,n-1\ldots,1)(0)}(x, y) = $\\
\hspace*{3cm}$
\prod\limits_{k=0}^{n-1}
[s_{n-1}\ldots s_1x, ~s_{0+k}d_ny]^{(-1)^k}[s_ns_{n-1}\ldots s_1(x),~s_{n}y].
$\\
Then 
$
[s_ns_{n-1}\ldots s_1(x),~s_{n}y]\equiv 1\qquad \mbox{mod}~(NG_n\cap D_n)d_{n+1}NG_{n+1}\cap D_{n+1}).
$
This gives the following if $\bar{x}\in C_1$ then $\bar{x}$ and $\partial_1\bar{x}$ acts on $C_n$ 
in the same way and so $\partial_1C_1$ acts trivially on $C_n.$\\
(iii) This axiom follows since 
$$
C_1({\bf G}) =  \dfrac{NG_1}{\partial_2(NG_2\cap D_2)}
=  \dfrac{NG_1}{[\mbox{Ker}d_1,~\mbox{Ker}d_0]}
$$
and $[\mbox{Ker}d_1, \mbox{Ker}d_0]$ contains the Peiffer elements so 
$(C_1({\bf G}), C_0({\bf G}), \partial)$ is a crossed module.

(iv) By defining
$$\partial_n\bar{z} = \overline{d_n^n(z)} \quad\mbox{with} ~~z\in NG_n$$
one obtains a well defined map $\partial: C_n({\bf G})\longrightarrow C_{n-1}({\bf G})$ satisfying
$\partial\partial =1.$ 
\end{pf}

One of the immediate consequences of the above is that if ${\bf {\mathbb
G}} = ({\mathbf G},f)$ is a simplicial group augmented over a group $G$, then
${\mathfrak{C}} = (C({\bf G}), f)$ is an augmented crossed complex over
$G$.  Moreover if ${\bf {\mathbb
G}}$ is exact at $G_0$, then ${\mathfrak{C}}$ is also exact at $C_0({\bf
G})$.  Thus to study what happens to a resolution we need only consider
the freeness and exactness in higher dimensions.

\section{`Step-by-Step' Constructions and CW-bases}
In this section, we describe the special case of the `step-by-step'
construction of a free simplicial resolution and its skeleton up to dimension
2 and will interpret this construction and see how that relates to other
constructions such as that of a free crossed module.

Many of the observations that we will make, do apply more generally to
arbitrary free simplicial groups with specified CW-basis, but our aim
here is limited to examining resolutions in some detail.
We first examine the relationship of a CW-basis to the step-by-step 
construction given earlier.

The 1-skeleton ${\bf{ F}}^{(1)}$ of a free simplicial resolution of
 a group $G$ was built by adding new indeterminates, for instance, 
in one to one 
correspondence with $\Omega^{1}$ a set of generators for $\pi_1({F}^{(0)})$,  
$F_1^1 = F_1^{(0)}(X_0) = F(s_0(X_0) \cup Y_1) \cong F(s_0(X_0))\ast F(Y_1)$, 
where $\ast$ is  free product, with the face maps and degeneracy map%
$$
\diagram
{F}(s_0(X_0)\cup Y_1)\rto<0.25ex>\rto<1ex>^{\qquad d_0,d_1 } & {F(X_0)}
\lto<0.75ex>^{\qquad s_0}\rto^{\quad d_0^0}&G   
\enddiagram
$$
where $F(X_0)\stackrel{d_{0}^{0}}{\longrightarrow} G$ is an augmentation map  and  $s_0,$ $d_0^1$ and $d_1^1$ 
are given by
$$
\begin{array}{ccc}
d_1^1(y_i)=b_i\in \text{Ker}d_0^0, & d_0^1(y_i)=1, & s_0(x_0)=s_0(x_0) ~~~\mbox{for}~~~ x_0\in X_0. 
\end{array}
$$
We note that this makes $\langle X_0 \mid d_1Y \rangle$ into a presentation of $G$ in the ordinary sense.
The 1-skeleton ${\bf {\bf{F}}}^{(1)}$ looks like:%
$$
\diagram{...\quad F(s_1s_0(X_0)\cup s_0(Y_1)\cup s_1(Y_1))
\rto<0.25ex> \rto<1ex> \rto<1.75ex>^{\qquad\qquad d_0 ,d_1 ,d_2 } & F(s_0(X_0)\cup (Y_1)) \lto<0.75ex> \lto<1.50ex>^{\qquad\qquad s_1,s_0}\rto<0.25ex> \rto<1ex>^{\qquad d_1, d_0} & F(X_0)\lto<0.75ex>^{\qquad s_0} .}\enddiagram
$$
Note that for $n>1,$ higher levels of ${\bf{F}}^{(1)}$ are generated by the
degenerate elements.
\begin{lem}\label{ber}
We assume given the 1-skeleton ${\bf{F}}^{(1)}$. Let $d_0^1$ and 
$d_{1\text{ }}^1$be evaluation homomorphisms. Then

i) \ \ {\rm Ker}$d_0^1= {\langle Y_1\rangle},$

ii) \ {\rm Ker}$d_1^1={\langle Z\rangle },$\\
where $Z=\{s_1(y)^{-1}s_0(y) : y \in Y_1\}$ and ${\langle Y_1\rangle}$
is normal closure of $Y_1.$
\end{lem}
\begin{pf} 
Clear.   
\end{pf}\\
Note $\pi_0({\bf{F} }^{(1)})\cong  G $.

The link between the bottom step of a step-by-step construction of the
resolution
 and that of a CW-basis $\mathfrak{F}$ is thus clear.
The 2-skeleton gives the non-degenerate elements of the resulting CW-basis,
 $\mathfrak{F}_2,$
and in general we can take $Y_n\cong \Omega^{n-1},$ and 
$\mathfrak{F}_n=Y_n\cup~\bigcup s_i(\mathfrak{F}_{n-1}).$ For both combinatorial and computational 
purposes, the way in which $Y_n$ corresponds to $\Omega^{n-1}$ can be important and in general 
it is necessary to specify the function
$
g^{n-1}:\Omega^{n-1}\longrightarrow NF_n
$
or its last face 
$
d_ng^{n-1}:\Omega^{n-1}\longrightarrow NF_{n-1}.
$

{\bf Remark:}
For homological and computational reasons, it is often useful also to specify the contracting homotopyon the underlying simplicial set of ${\mathbf F}$ and to build this into  the resolution progresses. We will not discuss how to do this here however as it is not needed for our immediate purposes. 

\medskip 

Before  carrying on the `step-by-step' construction of the free simplicial
group,   we will interpret the first homotopy group, $\pi _1({\bf {\mathbf F}}^{(1)})$,
of ${\mathbf F}^{(1)}$  to find what it looks like.

For  any simplicial group ${\mathbf F}$, if ${\mathbf F}={\mathbf F}^{(1)},$ 
then, 
$$ 
\pi_1({\mathbf F})={\rm Ker}(d_1 :{\rm Ker}d_0^1/[{\rm Ker}d_1^1, ~{\rm Ker}d_0^1] {\longrightarrow F_0}). 
$$
Indeed, by definition, the first homotopy group is  
$$
\pi_1({\mathbf F})=(\text{Ker}d_0^1\cap \text{Ker}d_1^1)/d_2^2(\text{Ker} d_0^2\cap \text{Ker}d_1^2). 
$$
By a lemma of Brown and Loday \cite{bl}, see also \cite{mp1}, the denominator  of this homotopy group isexactly  
$$
\partial_2(NF_2)=d_2^2(\text{Ker}d_0^2\cap \text{Ker}d_1^2)=[\text{Ker}d_0^1,  ~\text{Ker}d_1^1] 
$$
and the morphism  
$$ 
\delta :\text{Ker}d_0^1/\partial _2(NF_2)\longrightarrow F_0, 
$$
where  $\delta =d_1$ (restricted to $NF_1/\partial _2NF_2)$,  is a
crossed  module. Here $NF_0$ acts on $NF_1/\partial _2NF_2$ by conjugation via $%
s_0,$ that  is,
$$
\begin{array}{ccl}
NF_1/\partial_2NF_2\times NF_0 & \longrightarrow & NF_1/\partial _2NF_2, \\ 
({x},\overline{y})  & \longmapsto & {}^{x}
\overline{y}=\overline{s_0(x)ys_0(x)^{-1}},  
\end{array} 
$$
where  $\overline{y}$ denotes the corresponding element of $NF_1/\partial_2NF_2$ whilst $y\in NF_1.$ 
\begin{eqnarray*}
  \pi _1({\mathbf F}) & = &\text{Ker}(\text{Ker}d_0^1/\partial_2(NF_2)\longrightarrow F_0)\\
                    & = &\text{Ker}(NF_1/[\text{Ker}d_1, ~\text{Ker}d_0]\longrightarrow F_0). 
\end{eqnarray*}

\begin{prop} 
Given a presentation  $P=\langle{X_0~|~R\rangle}$ of a group $G$ and $
{\mathbf F}^{(1)}$, the 1-skeleton of the free simplicial group generated by
this presentation,  then
$$
\delta :NF_1^{(1)}/\partial_2(NF_2^{(1)})\longrightarrow NF_0^{(1)} 
$$
is  the free crossed module on $R\rightarrow {F(X_0)},$ and $\pi_2({\mathbf F}^{(1)})$
is the module of identities of the presentation $P.$
\end{prop}
\begin{pf}
Clear.
\end{pf}

\medskip

Note that for the case of ${\mathbf F}^{(2)}$, if $x_i,~x_j$ are in $NF_1^{(2)},$ then
generators of the normal subgroup $NF_2^{(2)}\cap D_2$ are of the form
$
[s_1(x_i)^{-1}s_0(x_j),~s_1(x_j)].
$
We now will recall the next step of the construction of a free simplicial
group. We take a set of generators 
$
\Omega ^1=\{S_i\}\subset \pi _1({\mathbf F}^{(1)}) 
$
and kill off the elements in the homotopy group $\pi _1({\mathbf F}^{(1)})$ by
adding new indeterminates $Y_2=\{y_i\}$ into $F_2^{(1)}$ where $Y_2$ is in $1-1$ 
correspondence with  $\Omega ^1$ to
establish 
$$
F_2^{(2)}=F_2^{(1)}(Y_2)=F(s_1s_0(X_0)\cup s_0(Y_1)\cup s_1(Y_1)\cup Y_2)  
$$
together with 
$$
d_0^2(y_i~)=1,\quad d_1^2(y_i~)=1,\quad d_2^2(y_i)=S_i,\textrm{ mod }\partial_3NF_3^{(2)}. 
$$
Hence the 2-skeleton ${\mathbf F}^{(2)}$ looks like 
$$
\diagram{ F(s_1s_0(X_0)\cup s_0(Y_1)\cup s_1(Y_1)\cup
  Y_2)\ar[r]<0.25ex> \ar[r]<1ex> \ar[r]<1.75ex>^{\qquad\qquad d_0 ,d_1 ,d_2 }&F(s_0(X_0)\cup (Y_1)) \ar[l]<0.75ex> \ar[l]<1.50ex>^{\qquad\qquad s_1,s_0}\ar[r]<0.25ex>\ar[r]<1ex>^{\qquad d_1, d_0} & F(X_0)\ar[l]<0.75ex>^{\qquad s_0} ,}\enddiagram
$$
and, of course,
$$
\begin{array}{lll}
F_2&=&F(s_1s_0(X_0)\cup s_0(Y_1)\cup s_1(Y_1)\cup Y_2)\\
&\cong&F(s_1s_0(X_0))\ast F(s_0(Y_1))\ast F(s_1(Y_1))\ast F(Y_2).
\end{array}
$$
In ${\bf F}^{(2)}$, higher levels than dimension 2 are generated by degenerate elements.

This pattern, of course, continues to higher dimensions. We thus have in 
each dimension, $k$-dimensional construction data $(\Omega^{k}, g^{(k)})$ and a 
$k^{th}$-level presentation of the group, $G.$ The various $\Omega^{k}$ thus provide us with a 
$CW$-basis for ${\bf F}.$
\section{Free Crossed Resolutions}
In this section we want to examine in slightly more detail this step-by-step construction
through the perspective of the corresponding crossed complex, examining not only to see
if $\mathfrak{C}(\mathbb{F})$ is a crossed resolution of a group $G,$ but also how the homotopy 
type of $\mathfrak{C}(\mathbb{F}^{(k)})$ is constructed from  $\mathfrak{C}(\mathbb{F}^{(k-1)}).$
Knowledge of this process would seem essential if the construction of crossed resolutions 
is to be `mechanised'. It also helps in the interpretation of homological invariants and 
their linkage with combinatorial properties of a presentation or of a higher level presentation 
of a group $G.$

As the analysis is applicable in greater generality, we start by looking at an arbitrary free 
simplicial group with chosen $CW$-basis.

A `step-by-step' construction of a free simplicial group is constructed from
simplicial group inclusions 
$$
{\bf F}^{(0)}\subseteq \text{ }{\bf F}^{(1)}\text{ }\subseteq 
\text{ }{\bf F}^{(2)}\text{ }\subseteq \text{ }\ldots 
$$
We take the functor $\mathfrak{C}$ which is described in
Section~\ref{fa}, to see what ${C}_n({\bf F}^{(k)})$ looks like, where 
${\bf  F}^{(k)}$ is the k-skeleton of that construction, concentrating
our attention in low dimensions.
For $k=0,$ we have the 0-skeleton ${\mathbb F}^{(0)}$ of the construction%
$$\diagram
{\mathbb F}^{(0)}:\ldots \quad {F(X_0)}\rto<.3ex>\rto<-.3ex>&
         {F(X_0)}\rightarrow {F(X_0)}/N. 
\enddiagram
$$
Here ${\bf F}^{(0)}$ is the trivial simplicial group in which in every
degree $n,$ $F_n^{(0)}= {F(X_0)}$ and $d_i^n=\ $id$\ =s_j^n.$
It is easy to see that ${C}_0({\bf{F}}^{(0)})= {F(X_0)}$ as $NF_1\cap D_1$ is 
trivial.
\noindent The 1-skeleton is 
$$
\diagram{
 ...{~}{F}(s_1s_0(X_0)\cup s_0(Y_1)\cup s_1(Y_1))
\rto<0.25ex> \rto<1ex> \rto<1.75ex>^{\qquad\qquad d_0 ,d_1 ,d_2 } 
& {F}(s_0(X_0)\cup Y_1) \lto<0.75ex> \lto<1.50ex>^{\qquad\qquad s_1,s_0} 
\rto<0.25ex> \rto<1ex>^{\qquad d_1, d_0} & {{F}(X_0)}
\lto<0.75ex>^{\qquad s_0},}\enddiagram
$$
and since $F_2^{(1)}$ is generated by  degenerate elements, $%
F_2^{(1)}=D_2$, so the crossed complex term ${C}_1({\bf   F}^{(1)})$ is the
following%
$$
\begin{array}{llll}
C_1({\bf F}^{(1)}) & = & \dfrac{NF_1^{(1)}}{(NF_1^{(1)}\cap D_1) \partial
_2(NF_2^{(1)}\cap D_2)}, &  \\  
\\
& = & \dfrac{NF_1^{(1)}}{\partial _2(NF_2^{(1)}\cap D_2)}\quad & \text{since 
}NF_1\cap D_1=1, \\  
\\
& = & \dfrac{NF_1^{(1)}}{\partial _2(NF_2^{(1)})}\quad
& \text{as }F_2^{(1)}=D_2. 
\end{array}
$$
By Lemma~\ref{ber} and the Brown-Loday lemma \cite{bl}, we have $NF_1^{(1)}=
\langle Y_1\rangle$ 
and $\partial _2(NF_2^{(1)})$ is generated by the Peiffer elements,
respectively. It then follows that 
$$
C_1({\bf F}^{(1)})=\langle Y_1\rangle /P_1.\ 
$$
Here $P_1$ is the first dimensional Peiffer normal subgroup. The proof of Theorem~\ref{rin}
from \cite{bh3} interprets within this context as showing that 
$
\partial _1:\langle Y_1\rangle/P_1\longrightarrow {F(X_0)} 
$
is the free crossed module on the presentation $\langle X_0\mid d_1(Y_1)\rangle$ of $\pi_0({\bf F}).$\\

Looking at the case $2,$ the 2-skeleton of the construction is%
$$
\diagram{
 ...{~}{F}(s_1s_0(X_0)\cup s_0(Y_1)\cup s_1(Y_1)\cup Y_2)
\rto<0.25ex> \rto<1ex> \rto<1.75ex>^{\qquad\qquad d_0 ,d_1 ,d_2 } 
& {F}(s_0(X_0)\cup Y_1) \lto<0.75ex> \lto<1.50ex>^{\qquad\qquad{~}~ s_1,s_0} 
\rto<0.25ex> \rto<1ex>^{\qquad d_1, d_0} & {{F}(X_0)}
\lto<0.75ex>^{\qquad s_0}.}\enddiagram
$$
As before $F_3^{(2)}=D_3$ as $F_3^{(2)}$ is generated by the degeneracy
elements. Thus the second term of the crossed complex is
$$
\begin{array}{lllc}
C_2({\bf  F}^{(2)}) & = & \dfrac{NF_2^{(2)}}{(NF_2^{(2)}\cap D_2)\partial
_3(NF_3^{(2)}\cap D_3)}, &  \\  
\\
& = & \dfrac{NF_2^{(2)}}{(NF_2^{(2)}\cap D_2)\partial _3(NF_3^{(2)})}\quad
& \text{as }F_3^{(2)}=D_3. 
\end{array}
$$
If $x,y\in NF_1,$ then $NF_2\cap D_2$ is generated by the elements of the
form
$
[s_1x^{-1}s_0x,~ s_1y] 
$
and in general, if $x,~y\in NF_{n-1},$ then
$
s_{n-1}x^{-1}s_{n-2}x,~ s_{n-1}y] \in NF_n\cap D_n. 
$
Now look at $\partial _3(NF_3^{(2)})$ in terms of the skeleton ${\bf  F}
^{(2)}.$ In a similar way to the proof of \thinspace Lemma~\ref{ber} and 
as $d_0^2(y_i)=d_1^2(y_i)=1 \ y\in Y_2$, \ one can readily  obtain that:
$$
NF_2^{(2)}= \langle s_1(Y_1)\cup Y_2\rangle\cap\langle Z\cup Y_2\rangle. 
$$
where $Z$ as in Lemma~\ref{ber}.

On the other hand, \cite{mp1} shows that on writing $K_I = \bigcap_{i\in I}
\text{Ker}d_i$ for $I \subseteq [n-1]$,
$$
\partial
_3(NG_3^{(2)})=\prod_{I,J} [K_I, ~K_J] [K_{\{0,2\}}, ~K_{\{0,1\}}] [K_{\{1,2\}}, ~K_{\{0,1\}}] 
[K_{\{1,2\}}, ~K_{\{0,2\}}], 
$$
where $I\cup J=[2],~I\cap J=\emptyset $, so this is generated by the following elements: for
$x_i\in NF_1=\text{Ker}d_0$ and $y_1,y_2\in NF_2=\,$Ker$d_0\cap
$Ker$d_1$ with $1\leq i,j\leq n$,
$$
\begin{array}{rcl}
\lbrack s_0{x_i}^{-1}s_1s_0d_1x_i,~ y_1 \rbrack &\quad(1)\\ 
\lbrack s_1{x_i}^{-1}s_0x_i,~ s_1d_2(y_1){y_1}^{-1} \rbrack &\quad(2) \\ 
\lbrack x_is_1d_2{x_i}^{-1}s_0d_2x_i,~ s_1y_1 \rbrack &\quad(3)\\
\lbrack {y_1}^{-1}s_1d_2y_1,~ y_2 \rbrack &\quad(4)\\ 
\lbrack y_1s_1d_2{y_1}^{-1}s_0d_2y_1,~ y_2 \rbrack &\quad(5) \\ 
\lbrack y_1s_1d_2{y_1}^{-1}s_0d_2y_1,~ s_1d_2(y_2){y_2}^{-1} \rbrack &\quad(6).
\end{array}
$$
The normal subgroup generated by these elements will be denoted by $P_2$ and will be
called the second dimensional Peiffer  normal subgroup. We thus in principle have not only an
explicit presentation of $C_2({\bf  F}^{(2)})$ but a list of seven `generic' 
moves analogous to the Peiffer moves introduced by Brown and Huebschmann, \cite{bh3}.

Writing $Q_2=NF_2^{(2)}\cap D_2,$ we get the second term of the crossed
complex as follows 
$$
C_2({\bf F}^{(2)})=\frac{\langle s_1(Y_1)\cup Y_2\rangle
\cap\langle Z\cup Y_2\rangle}{Q_2 \cdot P_2} 
$$
\begin{prop}
Let ${\bf F}^{(2)}$ be the 2-skeleton of a free simplicial group 
resolving $G= F(X_0)/N.$ Then 
$$
{\mathfrak{C}}^{(2)}:\quad NF_2^{(2)}
/(Q_2 \cdot P_2)\stackrel{\partial _2%
}{\rightarrow }\langle Y_1\rangle /P_1  
\stackrel{\partial _1}{\rightarrow }{F(X_0)}\stackrel{g}{
\rightarrow }{F(X_0)}/{N}\stackrel{f}{\rightarrow }1 
$$
is the 2-skeleton of a free crossed resolution of $G$ where $\partial _2$
and $\partial _1$ are given  respectively by: for $y_1\in\langle s_1(Y_1)\cup Y_2\rangle \cap \langle Z\cup Y_2\rangle$ and 
$x_i\in \langle Y_1\rangle,$
$$
\partial _2(y_1(Q_2\cdot P_2))=d_2(y_1) P_1\text{ and }\partial
_1(x_i P_1)=d_1(X_i). 
$$
where $NF_2^{(2)}$ is $(\langle s_1(Y_1)\cup Y_2\rangle\cap\langle Z\cup Y_2\rangle)$.
\end{prop}
\begin{pf}
This follows immediately from the description of the `step-by-step'
construction of the free simplicial group.
\end{pf}

\medskip

This result gives a combinatorial description of the $C_2({\bf F}^{(2)})$ term 
and if we manipulate the elements of $s_1(Y_1)\cup Y_2$ and $Z\cup Y_2$, remembering
that $Z=\{s_1(y)^{-1}s_0(y): y\in Y_1\},$ we can identify the generators as elements in the 
module of identities of the presentation $\langle X_0\mid d_1(Y_1)\rangle.$ The elements of $Y_2$ 
map via $d_2$ to a set of generators of this module since, of course, that is how they 
were chosen.

To complete our analysis of the r\^ole of a CW-basis in a free
simplicial resolution ${\mathbb F} = ({\bf F}, g)$ of a group $G$, we need
to check that $(C({\bf F}),C(g))$ is a free {\em crossed} resolution of
$G$ and to see what happens to the CW-basis in the `conversion'.

First a proposition showing how homotopies behave under the functor,
$\mathfrak C$.

\begin{prop}

Suppose ${\mathbf f_0}, {\mathbf f_1} :{\mathbf G}\rightarrow {\mathbf
H}$ are morphisms of simplicial groups and ${\mathbf h} : {\mathbf
f_0}\simeq{\mathbf f_1}$ is a homotopy between them.  Then ${\mathbf h}$
induces a homotopy ${\mathfrak C}({\mathbf h}) : \pi(1) \otimes{\mathfrak C}({\mathbf G})
\rightarrow{\mathfrak C}({\mathbf H})$ between ${\mathfrak C}{(\mathbf f_0})$ and ${\mathfrak C}({\mathbf
f_1})$.
\end{prop}
(Here the $\otimes $ is the tensor product of crossed complexes
introduced by Brown and Higgins, \cite{bh2}, and $\pi(1) :=
\pi(\Delta^1)$ is the groupoid `unit interval'.  For more on the
homotopy theory of crossed complexes, the simplicial category theory of
the category of crossed complexes, etc., see \cite{bgpt} and
\cite{tonks}.)

 \begin{pf}

The homotopy ${\mathbf h}$ can be realised as a morphism,
${\mathbf h} : \Delta[1] \bar{\otimes} {\mathbf G} \rightarrow {\mathbf
H},$ where $\Delta[1] \bar{\otimes} \quad $ is the simplicial tensor
within the simplicially enriched category of simplicial groups (or
groupoids) (see Quillen, \cite{q}, or the discussion in \cite{mp}.) This
is given as a colimit of copies of ${\bf G}$ by the construction outlined
in \cite{q}.

The functor ${\mathfrak C}$ can be thought of in two equivalent ways.  It is either
the composite  of the reflection onto the variety of simplicial
group(oid) T-complexes (cf. \cite{ep}) followed by the equivalence
between that and the category of crossed complexes, or alternatively it
uses the Cegarra-Carrasco equivalence between simplicial groupoids and
hypercrossed complexes of group(oid)s followed by the reflection onto
the variety of crossed complexes within that category. (The advantage at
this point in using groupoids is that $\pi(1)$ is naturally a groupoid,
but this can be avoided if desired.)  From either description it is
clear that ${\mathfrak C}$ will preserve colimits and thus tensors with simplicial
sets,
thus
$${\mathfrak C}(\Delta[1]\bar{\otimes} {\mathbf G}) \cong \Delta[1]\bar{\otimes}{\mathfrak C}({\mathbf G})\cong \pi(1)\otimes {\mathfrak C}({\mathbf G}).$$
Composing ${\mathfrak C}({\mathbf h})$ with these isomorphisms gives the
result.\end{pf}

\begin{cor}
If ${\mathbf g} : {\mathbf F} \rightarrow {\mathbf K}(G,0)$ is a free
simplicial resolution  of $G$, then ${\mathfrak C}({\mathbf g}) : {\mathfrak C}({\mathbf F}) \rightarrow {\mathfrak C}({\mathbf K}(G,0))= G$
is a free crossed resolution of $G$. 
\end{cor}
\begin{pf}
The data on ${\mathbf g}$ can be specified by giving a homotopy between
the identity on ${\mathbf F}$ and the map that `squashes $NF$ down to
$G$' and then uses a section of the augmentation map, $g_0$, to yield a map back to $NF_0$.    The corollary now
follows from the previous result applied to this simplicial homotopy.
\hfill\end{pf}

\medskip

To finish the comparison, we will show that each $C_n({\mathbf F})$ is a
free $G$-module on $Y_n$ if $n\geq 2$.  We start with $n = 2$ but in fact
almost the same proof works in higher dimensions.

Suppose that $M$ is a $G$-module and $\Theta :Y_2 \rightarrow M$ is a
function, we want to prove that $C_2({\mathbf F})$ is free on $Y_2$, so
we need to extend $\Theta$ to a map on $C_2({\mathbf F})$.  Form the
crossed complex
$$\ldots \rightarrow M \rightarrow 1\rightarrow G$$
with $M$ in dimension 2, $G$ in dimension 0, all other levels being
trivial and the action of $G$ on $M$ being the given one.  This has an associated 
simplicial group $S(M,G)$ with $NS(M,G)$ this crossed complex.    There
is an obvious morphism, $\phi$ from ${\mathbf F}^{(1)}$ to $S(M,G)$,
inducing the quotient morphism $g : F(X_0) \rightarrow F(X_0) / N \cong
G$.  As $\pi_1(S(M,G))$ is trivial, Proposition \ref{free} applies to
show $\phi$ extends over ${\mathbf F}^{(2)}$ also extending $\Theta$.  Now we 
use ${\mathfrak C}$ to pass back to crossed complexes to get $${\mathfrak C}(\phi) : {\mathfrak C}({\mathbf
F}^{(2)}) \rightarrow M$$ extending $\Theta$.  As $C_2({\mathbf F}^{(2)}) 
\simeq C_2({\mathbf F}) $, this proves the claim that $C_2({\mathbf
F}^{(2)}) $ is a free $G$-module on $Y_2$.

Of course, the only difference that is needed in dimension $n$ is in the
definition of $S(M,G)$, where $M$ is placed in dimension $n$ and $\Theta
: Y_n \rightarrow M$ is given.

We have proved:

\begin{prop}
If $\mathbb{F}$ is a simplicial resolution of $G$ given by a construction data sequence
$\{(Y_i, g^{(i)}), i=0,1,\ldots\}$ and ${\mathbf F}^{(k)}$ is the 
corresponding $k$-skeleton, 
then if $k\geq 2,$ $C_k(\mathbb{F}^{(k)})$ is a free $G$-module on $Y_k.$
\end{prop}

 Summarising we get:
\begin{thm}
The `step-by-step' construction of simplicial resolution of a group, $G$, yields a
`step-by-step' construction of a crossed resolution of $G$ via the crossed complex 
construction, $\mathfrak{C}.$
\end{thm}                   
As a bonus for our method we also have given an explicit description of
the crossed complex construction in low dimensions.  The construction
data to dimension $n$ yields an $n$-dimensional word system in the sense
of R.A.Brown.  What is less clear, as we have mentioned before, is why
the word system given by Whitehead (see \cite{rab}, Example 2.2.3) would not
seem to lift back to give construction data for a free simplicial group.

\noindent\begin{tabular}{lcl}
  A. MUTLU &\quad& T. PORTER \\
Department of Mathematics &\hspace{3mm}&School of Mathematics \\                   
Faculty of Science &&  Dean Street\\                               
University of Celal Bayar &&  University of Wales, Bangor\\               Muradiye Campus, 45030 && Gwynedd, LL57 1UT\\
Manisa, TURKEY &&  UK\\   
e-Mail: amutlu@spil.bayar.edu.tr  && e-Mail: t.porter@bangor.ac.uk\\
\end{tabular}


\begin{thebibliography}{99}
\bibitem{andre}
{\sc M. Andr\'{e},}  Homologie  des  Alg\`{e}bres  Commutatives,  
\ {\em  Lecture Notes in Math.,} {\em  Springer,} \ {\bf 206}, \
 \ (1970).

\bibitem{ashley}
{\sc N. Ashley,}
Simplicial T-Complexes: a non abelian version of a theorem of Dold-Kan,
Dissertationes Math., {\bf 165}, \ (1988), \ 11-58, \ 
{\em Ph.D. Thesis,}~  University of Wales, Bangor, \ (1978).

\bibitem{baues1}
{\sc H. J. Baues,} Algebraic Homotopy,
{\em Cambridge Studies in Advanced Mathematics}, {\bf 15},
Cambridge Univ. Press., (1989), 450 pages.

\bibitem{bgpt}
{\sc R. Brown}, {\sc M.Golasi\'nski}, {\sc T. Porter}, and {\sc A.Tonks,}
Spaces of maps into classifying spaces for equivariant crossed
complexes, {\em Indag. Mathem.}, N.S. {\bf 8(2)}, (1997), 157-172.

\bibitem{bh1}  
{\sc R. Brown} and {\sc P. J. Higgins,} 
Colimit-theorems for relative homotopy groups, 
{\em Jour. Pure Appl. Algebra,} \ {\bf 22}, \
(1981), \ 11-41.

\bibitem{bh2}  
{\sc R. Brown} and {\sc P. J. Higgins,} 
Homotopies and tensor products for $\omega$-groupoids and crossed complexes, 
{\em Jour. Pure Appl. Algebra,} \ {\bf 47}, \
(1987), \ 1-33.

\bibitem{bh3}
{\sc R. Brown} and {\sc J. Huesbschmann,}
Identities among relations, {\em Low dimension topology,}
London Math. Soc. Lecture Note Series, \ {\bf 48},
{\em (ed. R. Brown and T. L. Thickstun, Cambridge University Press)} 1982, pp. 153-202.

\bibitem{bl}  
{\sc R. Brown} and {\sc J.-L. Loday,}
Van Kampen Theorems for Diagrams of Spaces, {\em Topology},
\ {\bf 26}, \ (1987), \ 311-335.



\bibitem{rab}
{\sc R. A. Brown}, Generalized Group Presentations and Formal Deformations 
of CW complexes, {\em Trans. Amer. Math. Soc.}\ {\bf 334},\ (1992), \ 519-549.



\bibitem{carrasco}
{\sc P. Carrasco,}
Complejos    Hipercruzados,    Cohomologia    y 
Extensiones, {\em Ph.D. Thesis,} \ Universidad de Granada, \ (1987).



\bibitem{c:c}
{\sc P. Carrasco} and {\sc A. M. Cegarra,} Group-theoretic Algebraic Models 
for Homotopy Types, \ {\em Jour. Pure Appl. Algebra,}  \ {\bf 75}, \ (1991), 195-235.



\bibitem{conduche}
{\sc D. Conduch\'{e},}
Modules Crois\'{e}s G\'{e}n\'{e}ralis\'{e}s de Longueur 2,
{\em Jour. Pure Appl. Algebra,} \ {\bf 34}, \ (1984), \ 155-178.



\bibitem{curtis}
{\sc E. B. Curtis, }  Simplicial Homotopy Theory,\ 
{\em  Advances in Math.}, \  {\bf 6}, \ (1971), \ 107-209.



\bibitem{ep}
{\sc P. J. Ehlers} and {\sc T. Porter,}
Varieties of Simplicial Groupoids, I: Crossed Complexes.
{\em Jour. Pure Appl. Algebra,} {\bf 120}, (1997), 221-233; plus: Correction,
same journal (to appear).


\bibitem{kan2}
{\sc D. M. Kan,} 
A relation between CW-complexes and free c.s.s groups,
{\em Amer. Jour. Maths.,} ~ {\bf 81}, ~  (1959), ~ 512-528.



\bibitem{keune}
{\sc F. Keune,} Homotopical Algebra and Algebraic K-theory, {\em Thesis}, Universiteit  van Amsterdam, 1972.



\bibitem{may}
{\sc J. P. May,} Simplicial Objects in Algebraic Topology, \
{\em Van Nostrand,}  {\em Math. Studies}, \ {\bf 11}, \ (1967).



\bibitem{moore}
{\sc J. C. Moore,} Seminar in Algebraic Homotopy, {\em Princeton}, \ (1956).


\bibitem{mutlu}
{\sc A. Mutlu,}  Peiffer Pairings in the Moore Complex of a
Simplicial Group, {\em Ph.D. Thesis}, University of Wales, Bangor,
 \ (1997); Bangor Preprint, 97.11, available via
\textsf{http://www.bangor.ac.uk/ma/research/preprints/97prep.html} 



\bibitem{mp}
{\sc A. Mutlu} and {\sc T. Porter,}~  
Iterated Peiffer pairings in the Moore complex of a simplicial group, Applied
Categorical Structures (to appear);
 (1997), Bangor Preprint, 97.12, available via
\textsf{http://www.bangor.ac.uk/ma/research/preprints/97prep.html} 

\bibitem{mp1}
{\sc A. Mutlu} and {\sc T. Porter,}~  
Applications of Peiffer pairings in the Moore complex of a simplicial group, Theory and Applications of Categories,
4, No. 7, (1998) 148-173; previously as
 {\em Bangor Preprint} {\bf 97.17.} available via
 \textsf{http://www.bangor.ac.uk/ma/research/preprints/97prep.html} .


 \bibitem{mp2}{\sc A. Mutlu} and {\sc T. Porter,}~   Freeness Conditions for
   2-Crossed Modules and Complexes, Theory and Applications of Categories, 4, No.8, (1998) 174-194; previously as {\em Bangor Preprint} {\bf 97.19.} available via
\textsf{http://www.bangor.ac.uk/ma/research/preprints/97prep.html} .

 
\bibitem{porter}
{\sc T. Porter,}
 $n$-Types of simplicial groups and crossed $n$-cubes,
  {\em  Topology,} ~ {\bf 32}, ~ (1993), 5-24.



\bibitem{q}
{\sc D. Quillen,}
Homotopical Algebra,
\ {\em  Lecture Notes in Math.,} {\em  Springer,} \ {\bf 43}, \
 \ (1967).



\bibitem{tonks}
{\sc A. P. Tonks, }
Theory and Applications of crossed complexes, {\em Ph.D. Thesis},
 University of Wales, Bangor, (1993).


\bibitem{jh1}
{\sc J. H. C. Whitehead,}
Combinatorial Homotopy I and II,
{\em Bull. Amer. Math.  Soc.,} \ {\bf  55}, \  (1949), \  231-245 and 453-496.        

\end{thebibliography}
\end{document}